\documentclass[12pt,fleqn]{article}
\usepackage{amsmath}
\usepackage{amsfonts}
\usepackage{amssymb}

\newfont{\symb}{wncyr10 scaled\magstep1}

\newcommand{\C}{{\mathbb C}}

\newcommand{\R}{{\mathbb R}}

\normalbaselines
\newtheorem{theorem}{Theorem}[section]
\newtheorem{lemma}[theorem]{Lemma}

\newtheorem{definition}{Definition}[section]

\begin{document}

\title{$L^p$ solvability of the Stationary Stokes problem on domains with conical singularity in any dimension
}

\author{Martin Dindos\\
\small The University of Edinburgh and Maxwell Institute of Mathematical Sciences\\[-0.8ex]
\small JCMB, The King's Buildings Mayfield Road, Edinburgh EH9 3JZ, Scotland \\[-0.8ex]
\small \texttt{mdindos@ed.ac.uk}\\
\and
Vladimir Maz'ya\\
\small Department of Mathematical Sciences, M\&O Building,\\[-0.8ex]
\small University of Liverpool, Liverpool L69 7ZL, UK\\
\small and\\
\small Department of Mathematics,\\[-0.8ex]
\small Link\"oping University, SE-581 83 Link\"oping, Sweden\\[-0.8ex]
\small \texttt{vladimir.mazya@liu.se} }

\date{{In memory of Michael Sh. Birman}\\[0.8ex]
}

\abstract{The Dirichlet boundary value problem for the Stokes
operator with $L^p$ data in any dimension on domains with conical
singularity (not necessary a Lipschitz graph) is considered. We
establish the solvability of the problem for all $p\in
(2-\varepsilon,\infty]$ and also its solvability in
$C(\overline{D})$ for the data in $C(\partial D)$.}

\maketitle

\section{Introduction}

In this paper we study the Stokes system (which is the linearized
version of the stationary Navier-Stokes system) on a fixed domain
$D\subset{\R}^n$, for $n\ge 3$. In fact, we establish our result
for a both Lam\'e system ($\nu<1/2$) and the Stokes system
($\nu=1/2)$.

We want to consider a classical question of the solvability of the
$L^p$ Dirichlet problem on the domain $D$.

Let us recall that {\it the Dirichlet problem for the system
(\ref{eq01}) is $L^p$ solvable} on the domain $D$ if for all
vector fields $f\in L^p(\partial D)$ there is a pair of $(u,p)$
(here $u:D\to {\R}^n$, $p:D\to{\R}$) such that

\begin{eqnarray}\label{eq01}
-\Delta u+\nabla p&=&0,\qquad \mbox{div
}u+(1-2\nu)p=0\quad\mbox{in
}D,\\
u\big|_{\partial D}&=&f\quad\mbox{almost everywhere},\nonumber\\
u^*&\in& L^p(\partial D),\nonumber
\end{eqnarray}

and moreover for some $C>0$ independent of $f$ the estimate

$$\|u^*\|_{L^p(\partial D)}\le C\|f\|_{L^p(\partial D)}\qquad\mbox{holds.}$$

Here, the boundary values of $u$ are understood in the
nontangential sense, that is we take the limit
$$u\big|_{\partial D}(x)=\lim_{y\to x,\,y\in \Gamma(x)} u(y),$$
over a collection of interior nontangential cones $\Gamma(x)$ of
the same aperture and height and vertex at $x\in
\partial D$ and $u^*$ is the classical nontangential maximal function
defined as
$$u^*(x)=\sup_{y\in \Gamma(x)}|u(y)|,\qquad\mbox{for all }x\in\partial D.$$

Furthermore, we say that the Dirichlet problem (\ref{eq01}) is
solvable for continuous data, if for all $f\in C(\partial D)$ the
vector field $u$ belongs to $C(\overline{D})$ and the estimate
$$\|u\|_{C(\overline{D})}\le C\|f\|_{C(\partial D)}\qquad\mbox{holds.}$$

So far, we have not said anything about the domain $D$, except the
requirement on the existence of the interior nontangential cones
$\Gamma(x)\subset D$ of the same aperture and height at every
boundary point $x\in \partial D$. For symmetry reasons let us also
assume the existence of exterior (i.e. in
${\R}^n\setminus\overline{D}$) nontangential cones, as well as
that the domain $D$ is bounded.

A most classical example of a domain $D$ that satisfies these
assumptions is a Lipschitz domain which is a domain for which the
boundary $\partial D$ can be locally described as a graph of a
Lipschitz function.

Another important class of domains satisfying outlined assumptions
are so-called generalized polyhedral domains (see \cite{MP2} and
\cite{MR}). The precise definition is rather complicated and
unnecessary for our purposes what we have in mind are domains that
look like polyhedra, however we allow the sides or edges to be
curved not just flat. For example, in two dimensions the boundary
of such domain will consist of a finite set of vertices joined by
$C^1$ curves meeting at the vertices nontransversally.

At the first sight, one might assume that the class of generalized
polyhedral domains is a subset of the class of Lipschitz domains.
This is however only true when $n=2$, when $n>3$ it is no longer
the case. Is it however clear that these classes are related.

The $L^p$ Dirichlet boundary value problem on Lipschitz domains
for the Laplacian has a long history starting with pioneering work
of Dahlberg \cite{D}, Jerison and Kenig \cite{JK}. As follows from
these results the $L^p$ Dirichlet boundary problem for the
Laplacian is solvable for all $p\in (2-\varepsilon,\infty]$
regardless of dimension. There are two key ingredients to
establish this result: the so called Rellich estimates (when
$p=2$) and the maximum principle (when $p=\infty$). Interpolating
these two leads to the stated result.

The Stokes (and Lam\'e) equations are PDE systems not a single
(scalar) equation. This implies that the second ingredient of the
above approach is not readily available since the maximum
principle that holds for PDE equations is not applicable to
general PDE systems. In addition, if we go outside the class of
Lipschitz domains, the Rellich estimates are also not available.

Despite that in low dimensions $n=2,3$ a weak version of the
maximum principle does hold \cite{Sh} which still allows to prove
$L^p$ solvability for all $p\in (2-\varepsilon,\infty]$, provided
the domain is Lipschitz. See \cite{DK2} or \cite{DM} for most
up-to date approach that works even on Riemannian manifolds. See
also \cite{BS} for regularity issues related to the Stokes system
in Lipschitz domains.

The question whether the same is true if $n>3$ is open and only
partial range of $p$ for which the problem is solvable is known.

In this paper we consider the problem outlined above for domains
with a single singular point, a conical vertex. Apart from this
point our domain will be smooth. In three dimensions the situation
has been considered in a different context before \cite{MP} where
estimates such as (\ref{eq15}) in three dimensions are
established. We mention also the recent book \cite{MR} where
strong solutions of the three-dimensional problem (\ref{eq01}) and
stationary Navier-Stokes equation are studied in detail. Our
approach is general and works in any dimension $n$. In principle
what we present here is fully extendable to all generalized
polyhedral domains. However to avoid technical challenges we only
focus on the particular case of an isolated conical singularity.

Our main result shows that in the setting described above (a
domain with one conical point)  the range of solvability $p\in
(2-\varepsilon,\infty]$ remains true in {\it any} dimension.

As we stated above the classes of Lipschitz and polyhedral domains
are related but neither is a subset of the other. However the
result we present is a strong indication that the range $p\in
(2-\varepsilon,\infty]$ should hold even for Lipschitz domains. In
fact, the known counterexamples to solvability in $L^p$ are shared
by these two classes of domains (see \cite{K} for such examples
when $p<2$).

The paper is organized as follows. In section 2 we establish
estimates for eigenvalues of a certain operator pencil for Lam\'e
and Stokes systems in a cone that holds in any dimension. These
estimates are in the spirit of work done in \cite{MP}. In section
3 we prove estimates for Green's function that are based upon
section 2 and finally in section 4 we present our main result
Theorem \ref{t2} and its proof that is based upon the explicit
estimates for the Green's function from section 3 and
interpolation.

\section{The operator pencil for Lam\'e and Stokes systems in a cone}

Let $\Omega\subset {\mathbb S}^{n-1}$. Suppose that cap$({\mathbb
S}^{n-1}\setminus\Omega)>0$. Then the first eigenvalue of the
Dirichlet problem for the operator $-\Delta_{{\mathbb S}^{n-1}}$
in $\Omega$ is positive. We represent this eigenvalue in the form
$M(M+n-2)$, $M>0$.

Consider the cone ${\cal K}=\{x\in {\R}^n; 0<|x|<\infty,
x/|x|\in\Omega\}.$ Or goal is to understand the solutions of the
system (\ref{eq1a})-(\ref{eq2a})
\begin{equation}\label{eq1a}
-\Delta U+\nabla P=0,\qquad \mbox{div }U+(1-2\nu)P=0\quad\mbox{in
}{\cal K}
\end{equation}
with the boundary condition
\begin{equation}\label{eq2a}
U\big|_{\partial K\setminus\{0\}}=0,
\end{equation}
in the form $U(x)=r^{\lambda_0}u(\omega),\,
P(x)=r^{\lambda_0-1}p(\omega)$. This requires study of spectral
properties of a certain operator pencil ${\cal L(\lambda)}$
defined below (\ref{eqPP}).

Here $\nu$ is the so-called the Poisson ratio. If $\nu<1/2$ the
equation (\ref{eq1a}) can be written in more classical
(elasticity) form
\begin{equation}\label{eq1b}
\Delta U+(1-2\nu)^{-1}\nabla\nabla\cdot U=0\quad\mbox{in }{\cal
K}.
\end{equation}

The case $\nu=1/2$ corresponds to the Stokes system
\begin{equation}\label{eq1}
-\Delta U+\nabla P=0,\qquad \mbox{div }U=0\quad\mbox{in }{\cal K}
\end{equation}
with the boundary condition
\begin{equation}\label{eq2}
U\big|_{\partial K\setminus\{0\}}=0.
\end{equation}

Now we define the pencil. We write (\ref{eq1a}) in the polar form
for $U=r^{\lambda}(u_u,u_\omega)$. After multiplying by
$r^{2-\lambda}$ we obtain:
\begin{eqnarray}\nonumber
-\Delta_{{\mathbb
S}^{n-1}}u_r&-&(\lambda+1)(\lambda+n-1)u_r-\frac{\lambda-1}{1-2\nu}[(\lambda+n-1)u_r\\&+&\nabla_{\omega}\cdot
u_\omega]+2[(\lambda+n-1)u_r+\nabla_{\omega}\cdot
u_\omega]=0\label{eqP},
\end{eqnarray}
and
\begin{eqnarray}\nonumber
Lu_\omega&-&(\lambda+1)(\lambda+n-1)u_\omega-\frac1{1-2\nu}[(\lambda+n-1)\nabla_\omega
u_r\\&+&\nabla_\omega(\nabla_\omega\cdot
u_\omega)]+2[(\lambda+n-1)u_\omega-\nabla_{\omega}
u_r]=0\label{eqP2},
\end{eqnarray}
where $L$ is the second order differential operator acting on
vector fields on ${\mathbb S}^{n-1}$ that arises from the
Navier-Stokes equation on the sphere ($L=-\Delta_{{\mathbb
S}^{n-1}}-\nabla_{\omega}(\nabla_\omega\cdot\,)-2$ and
$\Delta_{{\mathbb S}^{n-1}}$ is the Hodge Laplacian on ${\mathbb
S}^{n-1}$) (see also section 3.2.3 of \cite{KMR} for the
corresponding calculation in three dimensions).

Hence we arrive at the matrix differential operator ${\cal
L}(\lambda)\left(\begin{array}{c}u_r\\u_\omega
\end{array}\right)=0$, where:
\begin{eqnarray}
{\cal L}(\lambda)\left(\begin{array}{c}u_r\\u_\omega
\end{array}\right)=\left(\begin{array}{c}-\Delta_{{\mathbb
S}^{n-1}}u_r-\frac{2-2\nu}{1-2\nu}(\lambda-1)(\lambda+n-1)u_r+\frac{3-4\nu-\lambda}{1-2\nu}\nabla_{\omega}\cdot u_\omega\\
L_\nu
u_\omega-(\lambda-1)(\lambda+n-1)u_\omega-\frac{n+1-4\nu+\lambda}{1-2\nu}\nabla_\omega
u_r
\end{array}\right).\nonumber\\\mbox{ }\label{eqPP}
\end{eqnarray}

Here $L_\nu=L-(1-2\nu)^{-1}\nabla_\omega(\nabla_\omega\cdot\,\,)$.
This defines the operator pencil. The pencil operator is
self-adjoint if and only if Re $\lambda=-\frac{n-2}{2}$.

From now on we consider arbitrary $\nu\le 1/2$. In particular,
this includes the case of the Stokes system we are mostly
interested in.

For arbitrary real $t$ consider
\begin{equation}\label{phi}
\phi(t)=(t+1)(t+n-1)(2t+n-2)-(3-4\nu-t)(M-t)(M+t+n-2).
\end{equation}
Let $t(M)$ be the smallest solution of the equation
$$\phi(t)=(n-1)(2t+n-2)$$
in the interval $-(n-2)/2<t<M$. We claim that $t(M)>0$.

Indeed, for $t$ in the interval $[-(n-2)/2,0]$ we have
$$\phi(t)<(t+1)(t+n-1)(2t+n-2)<(n-1)(2t+n-2),$$
hence there is no solution in this interval. On the other hand
$\phi(M)=(M+1)(M+n-1)(2M+n-2)>(n-1)(2M+n-2)$, hence the solution
always exists and is positive.

The following result extends a similar three-dimensional result in
\cite{MP} (see also Section 5.5.4 in \cite{KMR}).

\begin{theorem}\label{t1} The strip determined by the inequality
\begin{equation}
\left|\mbox{Re }\lambda+\frac{n-2}2\right|<
\min\{1,t(M)\}+\frac{n-2}2\label{eq3}
\end{equation}
does not contain any eigenvalues of the pencil ${\cal L}(\lambda)$
in $\cal K$, provided $\nu\le1/2$. What this means is that the
boundary value problem (\ref{eq1a})-(\ref{eq2a}) has no solution
of the form
\begin{equation}
U(x)=r^{\lambda_0}u(\omega),\qquad
P(x)=r^{\lambda_0-1}p(\omega),\label{eq4}
\end{equation}
for $\lambda_0$ in this strip.
\end{theorem}

\noindent{\it Proof.} We will only consider the case when Re
$\lambda>-(n-2)/2$. This is enough, since by Theorem 3.2.1 of
\cite{KMR} $\lambda_0$ is an eigenvalue of the pencil if and only
if $-(n-2)-\overline{\lambda_0}$ is an eigenvalue. Also there is
no eigenvalue on the line Re $\lambda=-(n-2)/2$.

Let us therefore consider a pair $(U,P)$ of the form (\ref{eq4})
that solves the boundary value problem (\ref{eq1a})-(\ref{eq2a})
in $\cal K$. For arbitrary $\varepsilon>0$ we consider the domain
$${\cal K}_\varepsilon=\{x\in{\cal K}:\, \varepsilon<|x|<1/\varepsilon\}.$$
Then (\ref{eq1}) implies
$$-\int_{{\cal K}_\varepsilon}\Delta U\cdot{\overline U} dx+\int_{{\cal K}_\varepsilon}\nabla P\cdot{\overline U}dx=0.$$
Integrating by parts and using the equation for the divergence of
$U$ we obtain that
\begin{eqnarray}
&&\int_{{\cal K}_\varepsilon}|\nabla U|^2
dx+\varepsilon^{n-1}\int_{r=\varepsilon}{\overline
U}\cdot\partial_rU
d\omega-\varepsilon^{-n+1}\int_{r=1/\varepsilon}{\overline
U}\cdot\partial_rU d\omega\\&& +(1-2\nu)\int_{{\cal
K}_\varepsilon}|P|^2dx-\varepsilon^{n-1}\int_{r=\varepsilon}P{\overline
U_r} d\omega+\varepsilon^{-n+1}\int_{r=1/\varepsilon}P{\overline
U_r} d\omega=0.\nonumber
\end{eqnarray}
Here $U_r$ denotes the radial part of the vector $U$ and $d\omega$
is the surface measure on $S^{n-1}$. Now using the assumption
(\ref{eq4}) we obtain
\begin{eqnarray}\label{eq5}
&&\int_{\varepsilon}^{1/\varepsilon}r^{2Re\,\lambda_0+n-3}dr\int_{\Omega}\left(|\nabla_{\omega}u|^2+|\lambda_0|^2|u|^2+(1-2\nu)|p|^2\right)d\omega\\
&&+\left(\varepsilon^{2Re\,\lambda_0+n-2}-\varepsilon^{-(2Re\,\lambda_0+n-2)}\right)\left(\lambda_0\int_{\Omega}|u|^2d\omega-
\int_{\Omega}p{\overline u_r}d\omega\right)=0.\nonumber
\end{eqnarray}
We take the real part of equation (\ref{eq5}) and then integrate
in $r$. This gives us
\begin{eqnarray}\nonumber
\int_{\Omega}|\nabla_{\omega}u|^2d\omega &+&
\left((Re\,\lambda_0)^2+(Im\,\lambda_0)^2-(2Re\,\lambda_0+n-2)Re\,\lambda_0
\right)\int_{\Omega}|u|^2d\omega\\+(1-2\nu)\int_{\Omega}|p|^2dx&+&(2Re\,\lambda_0+n-2)Re\int_{\Omega}p{\overline
u_r}d\omega=0,\nonumber
\end{eqnarray}
which can be simplified to
\begin{eqnarray}\nonumber
\int_{\Omega}|\nabla_{\omega}u|^2d\omega &+&
\left((Im\,\lambda_0)^2-Re\,\lambda_0(Re\,\lambda_0+n-2)
\right)\int_{\Omega}|u|^2d\omega\\+(1-2\nu)\int_{\Omega}|p|^2dx&+&(2Re\,\lambda_0+n-2)Re\int_{\Omega}p{\overline
u_r}d\omega=0.\label{eq7}
\end{eqnarray}

Now we use the original equations to get an expression for the
pressure $p$. It follows that
$$r\partial_rP=x\cdot\nabla P=x\cdot\Delta U=\Delta(x\cdot U)-2\mbox{ div } U=\Delta
(rU_r)+2(1-2\nu)P.$$ Again using (\ref{eq4}) we obtain that

\begin{eqnarray}\nonumber
r\partial_rP&=&r^{\lambda_0-1}(\lambda_0-1)p(\omega),\\
\Delta(rU_r)&=&\Delta(r^{\lambda_0+1}u_r(\omega))=
r^{\lambda_0-1}\left(\Delta_{{\mathbb
S}^{n-1}}u_r(\omega)+(\lambda_0+1)(\lambda_0+n-1)u_r(\omega)\right).\nonumber
\end{eqnarray}

Hence $p(\omega)=-(3-4\nu-\lambda_0)^{-1}(\Delta_{{\mathbb
S}^{n-1}}u_r(\omega)+(\lambda_0+1)(\lambda_0+n-1)u_r(\omega))$.
Consequently,
$$Re\int_{\Omega}p{\overline
u_r}\,d\omega=Re\,(3-4\nu-\lambda_0)^{-1}\int_{\Omega}|\nabla_\omega
u_r|^2d\omega-Re\frac{(\lambda_0+1)(\lambda_0+n-1)}{3-4\nu-\lambda_0}\int_{\Omega}|u_r|^2d\omega.$$
From this identity and (\ref{eq7}) after multiplying by
$|3-4\nu-\lambda_0|^2$ we obtain:

\begin{eqnarray}\nonumber
0&=&|3-4\nu-\lambda_0|^2\left[\int_{\Omega}|\nabla_{\omega}u|^2d\omega
+ \left((Im\,\lambda_0)^2-Re\,\lambda_0(Re\,\lambda_0+n-2)
\right)\int_{\Omega}|u|^2d\omega\right]\\
&+&|3-4\nu-\lambda_0|^2(1-2\nu)\int_{\Omega}|p|^2dx\nonumber\\
&+&(2Re\,\lambda_0+n-2)(3-4\nu-Re\,\lambda_0)\int_{\Omega}|\nabla_\omega
u_r|^2d\omega\label{eq8}\\&-&(2Re\,\lambda_0+n-2)Re\left[(\lambda_0+1)(\lambda_0+n-1)(3-4\nu-{\overline\lambda_0})\right]\int_{\Omega}|u_r|^2d\omega.\nonumber
\end{eqnarray}

Now we use the fact that $M(M+n-2)$ is the first eigenvalue of the
Dirichlet problem for Laplacian on $\Omega$. Hence the inequality
$$M(M+n-2)\int_{\Omega}|u|^2d\omega\le \int_{\Omega}|\nabla_\omega u|^2d\omega$$
holds. Also
\begin{eqnarray}\nonumber
&&Re\left[(\lambda_0+1)(\lambda_0+n-1)(3-4\nu-{\overline\lambda_0})\right]\\&=&(Re\,\lambda_0+1)(Re\,\lambda_0+n-1)(3-4\nu-{Re\,\lambda_0})-|Im\,\lambda_0|^2(Re\,\lambda_0+n+3-4\nu)
\nonumber\\
&\le&(Re\,\lambda_0+1)(Re\,\lambda_0+n-1)(3-4\nu-{Re\,\lambda_0}).\nonumber
\end{eqnarray}
Using these two inequalities and the fact that
$2Re\,\lambda_0+n-2\ge 0$ we get that
\begin{eqnarray}\nonumber
&&|3-4\nu-\lambda_0|^2\left[M(M+n-2)
-Re\,\lambda_0(Re\,\lambda_0+n-2)
\right]\int_{\Omega}|u|^2d\omega\\
&+&|3-4\nu-\lambda_0|^2(1-2\nu)\int_{\Omega}|p|^2dx\nonumber\\
&+&(2Re\,\lambda_0+n-2)(3-4\nu-Re\,\lambda_0)\int_{\Omega}|\nabla_\omega
u_r|^2d\omega\label{eq9}\\&\le&(2Re\,\lambda_0+n-2)(Re\,\lambda_0+1)(Re\,\lambda_0+n-1)(3-4\nu-{Re\,\lambda_0})
 \int_{\Omega}|u_r|^2d\omega.\nonumber
\end{eqnarray}
We obtain further simplification by using the inequality
$|3-4\nu-\lambda_0|^2\ge (3-4\nu-Re\,\lambda_0)^2$. This gives us
\begin{eqnarray}\nonumber
&&(3-4\nu-Re\,\lambda_0)(M-Re\,\lambda_0)(M+Re\,\lambda_0+n-2)\int_{\Omega}|u|^2d\omega\\\label{eq10}
&+&(3-4\nu-Re\,\lambda_0)(1-2\nu)\int_{\Omega}|p|^2dx\nonumber\\
&+&(2Re\,\lambda_0+n-2)\int_{\Omega}|\nabla_\omega
u_r|^2d\omega\\&\le&(2Re\,\lambda_0+n-2)(Re\,\lambda_0+1)(Re\,\lambda_0+n-1)
 \int_{\Omega}|u_r|^2d\omega.\nonumber
\end{eqnarray}

We write $|u|^2=|u_r|^2+|u_\omega|^2$. This gives

\begin{eqnarray}\nonumber
&&(3-4\nu-Re\,\lambda_0)(M-Re\,\lambda_0)(M+Re\,\lambda_0+n-2)\int_{\Omega}|u_\omega|^2d\omega\\
&+&(3-4\nu-Re\,\lambda_0)(1-2\nu)\int_{\Omega}|p|^2dx\nonumber\\
&+&(2Re\,\lambda_0+n-2)\int_{\Omega}|\nabla_\omega
u_r|^2d\omega\le\phi(Re\,\lambda_0)\int_{\Omega}|u_r|^2d\omega,\label{eq11}\end{eqnarray}
where $\phi$ is defined by (\ref{phi}).

Using the equation div $U+(1-2\nu)P=0$ we get that
$\nabla_\omega\cdot u_\omega+(\lambda_0+n-1)u_r+(1-2\nu)p=0$.
Integrating this over $\Omega$ we conclude that
\begin{equation}
\int_{\Omega}\left(u_r+\frac{1-2\nu}{\lambda_0+n-1}p\right)d\omega=0.
\end{equation}

Consider therefore the following minimization problem for the
functional
\begin{equation}
\int_{\Omega}|\nabla_\omega
v|^2\,d\omega+\frac{(3-4\nu-Re\,\lambda_0)(1-2\nu)}{2Re\,\lambda_0+n-2}\int_{\Omega}|q|^2d\omega\label{min}
\end{equation}
for pairs of functions $(v,q)\in W^{1,0}_2(\Omega)\times
L^2(\Omega)$ such that
\begin{equation}
\int_{\Omega}\left(v+\frac{1-2\nu}{\lambda_0+n-1}q\right)d\omega=0\quad\mbox{
and}\quad\|v\|_{L^2(\Omega)}=1.\label{restr}
\end{equation}

 Let us denote by $\Theta(\Omega,\lambda)$
that minimum of this functional and let $(v_0,q_0)$ be a pair of
functions realizing this minimum. We choose $q_1$ to be arbitrary
$L^2(\Omega)$ function orthogonal to $1$. Then the pair
$(v_0,q_0+\alpha q_1)$ for any $\alpha\in\R$ satisfies
(\ref{restr}). Inserting this pair into (\ref{min}) we obtain
$$Re \int_{\Omega} q_0\cdot\overline{q_1}\,d\omega=0.$$
Consequently $q_0$ is a constant and we can restrict ourselves to
constant $q$ in the above formulated variational problem. From
this
\begin{eqnarray}
\Theta(\Omega,\lambda)&=&\inf\left\{\int_{\Omega}|\nabla_\omega
v|^2\,d\omega+\frac{(3-4\nu-Re\,\lambda_0)|\lambda_0+n-1|^2}{(1-2\nu)(2Re\,\lambda_0+n-2)|\Omega|}\left|\int_{\Omega}v\,d\omega\right|^2
\right\},\nonumber\\&&\label{th1}
\end{eqnarray}
where the infimum is taken over all
$W_2^{1,0}(\Omega)$ functions $v$ with $L^2(\Omega)$ norm $1$.

When $\nu=1/2$ this minimization problem takes slightly different
form and simplifies to
\begin{equation}
\Theta(\Omega)=\inf\left\{\int_{\Omega}|\nabla_\omega
v|^2\,d\omega:\,\int_{\Omega}v\,d\omega=0,\,v\big|_{\partial\Omega}=0\mbox{
and }\|v\|_{L^2(\Omega)}=1\right\}.\label{th2}
\end{equation}
Note that is this case the minimization problem is independent of
$\lambda$.

Considering Re $\lambda_0\in\left(-\frac{n-2}{2},3-4\nu\right]$ we
see the number
$$C=\frac{(3-4\nu-Re\,\lambda_0)|\lambda_0+n-1|^2}{(1-2\nu)(2Re\,\lambda_0+n-2)|\Omega|}>0.$$
Writing $v$ for the problem (\ref{th1}) in the form $v=c_0+v_1$,
where $c_0$ is a constant function and $v_1$ has average $0$ over
$\Omega$ we see that
$$\int_{\Omega}|\nabla_\omega
v|^2\,d\omega+C\left|\int_{\Omega}v\,d\omega\right|^2=\int_{\Omega}|\nabla_\omega
v_1|^2\,d\omega+C|\Omega|^2c_0^2.$$ Using (\ref{th2}) the integral
$\int_{\Omega}|\nabla_\omega v_1|^2\,d\omega$ can be further
estimated from below by
$\Theta(\Omega)\|v_1\|^2_{L^2(\Omega)}=\Theta(\Omega)(1-\|c_0\|^2_{L^2(\Omega)})$.

Hence
$$\Theta(\Omega,\lambda)=\inf\{\Theta(\Omega)(1-\|c_0\|^2_{L^2(\Omega)})+C|\Omega|^2c_0^2;\,c_0\in\C\mbox{ and
}\|c_0\|^2_{L^2(\Omega)}\le 1\}.$$ It follows that the infimum
will be attained either when $c_0=0$ or $\|c_0\|_{L^2(\Omega)}=1$.
Hence
$$\Theta(\Omega,\lambda)=\min\{\Theta(\Omega),C|\Omega|\},$$
or
\begin{equation}
\Theta(\Omega,\lambda)=\min\left\{\Theta(\Omega),\frac{(3-4\nu-Re\,\lambda_0)|\lambda_0+n-1|^2}{(1-2\nu)(2Re\,\lambda_0+n-2)}\right\}.
\label{th3}\end{equation}

Looking at (\ref{th2}) we can see that the function $\Omega\mapsto
\Theta(\Omega)$ does not increase as $\Omega$ increases. It
follows that $\Theta(\Omega)\ge\Theta({\mathbb S}^{n-1})$. This
minimum for the $n-1$ dimensional sphere ${\mathbb S}^{n-1}$ is
known explicitly and equals to $n-1$. This and (\ref{th3}) implies
that for Re $\lambda\le 1$ we have that $\Theta(\Omega,\lambda)\ge
n-1$.

Using this for $v=u_r$ in (\ref{eq11}) we see that
\begin{eqnarray}\nonumber
&&(3-4\nu-Re\,\lambda_0)(M-Re\,\lambda_0)(M+Re\,\lambda_0+n-2)\int_{\Omega}|u_\omega|^2d\omega\\&\le&
[\phi(Re\,\lambda_0)-(n-1)(2
Re\,\lambda_0+n-2)]\int_{\Omega}|u_r|^2d\omega.\label{eq12}\end{eqnarray}

It follows that for all $-(n-2)/2<\mbox{Re}\,\lambda_0<t(M)$ the
term on the righthand side of the inequality (\ref{eq12}) is not
positive, but the term on the lefthand side is nonnegative. Hence
both terms have to vanish, i.e., $u_r=u_\omega=0$ or $u$ is
constant. Given that $u$ vanishes on $\partial\Omega$ we get that
$U(x)=r^{\lambda_0}u(\omega)=0$ everywhere. This establishes the
claim.

\begin{theorem}\label{t3} Consider any (energy) solution of the
system (\ref{eq01}) in ${\cal K} \cap B(0,1)$ such that
$$u\big|_{\partial {\cal K}\cap B(0,1)}=0,\qquad u\big|_{\partial B(0,1)\cap {\cal K}}\in C(\Omega).$$
Then
\begin{equation}
u\in C^{\alpha}(\overline{{\cal K}\cap
B(0,1/2)}),\quad\mbox{and}\quad |\nabla u(x)|\le
C|x|^{\alpha-1}\|u\|_{C(\overline{\partial B(0,1)\cap {\cal
K}})},\label{eq335}
\end{equation}
for all $|x|\le 1/2$ and some $\alpha>0$, $C>0$ independent of
$u$.

Similarly, if $u$ is an (energy) solution of the system
(\ref{eq01}) in ${\cal K} \setminus {\overline {B(0,1)}}$ such
that
$$u\big|_{\partial {\cal K}\setminus B(0,1)}=0,\qquad u\big|_{\partial B(0,1)\cap {\cal K}}\in C(\Omega).$$
Then
\begin{equation}|u(x)|\le C|x|^{2-n-\alpha}\|u\|_{C(\overline{\partial
B(0,1)\cap {\cal K}})},\quad |\nabla u(x)|\le
C|x|^{1-n-\alpha}\|u\|_{C(\overline{\partial B(0,1)\cap {\cal
K}})},\label{eq336}
\end{equation}
for all $|x|\ge 2$ and some $\alpha>0$, $C>0$ independent of $u$.
\end{theorem}

\noindent{\it Proof.} This is a consequence of Theorem \ref{t1}.
Indeed, we first look for solutions of (\ref{eq01}) in the
separated form $U(x)=r^{\lambda}u(\omega)$,
$P(x)=r^{\lambda-1}p(\omega)$. Choose, any $\alpha\in
(0,\min\{1,t(M)\})$. Then by Theorem \ref{t1}, Re $\lambda>\alpha$
and any (energy) solution of the system (\ref{eq01}) on ${\cal
K}\cap B(0,1)$ has the following asymptotic representation by the
Kondrat'ev's theorem \cite{Ko}:
\begin{equation}
u(r,\omega)=
r^{\lambda_0}\sum_{i=1}^m\sum_{j=0}^{k_i}c_{ij}\frac1{j!}(\log
r)^ju^i_{j}(\omega) +O(r^{\text{Re }
\lambda_0-\varepsilon}),\qquad\mbox{for all
}\varepsilon>0.\label{eq334}
\end{equation}
Here $u^i_0(\omega)$ are the eigenfunctions of the pencil ${\cal
L}(\lambda_0)$ corresponding to an eigenvalue $\lambda_0$ with the
smallest real part ($m$ is the multiplicity of this eigenvalue)
and $u^i_1(\omega),u^i_2(\omega),\dots,u^i_{k_i}(\omega)$ are the
generalized eigenfunctions (a Jordan chain). Given our assumption
that $\partial\Omega$ is smooth, these are $C^\infty$
vector-valued functions.

Similarly, in the second case working on ${\cal K} \setminus
{\overline {B(0,1)}}$ we have
\begin{equation}
u(r,\omega)=
r^{2-n-\lambda_0}\sum_{i=1}^m\sum_{j=0}^{k_i}c_{ij}\frac1{j!}(\log
r)^ju^i_{j}(\omega)+O(r^{2-n-\text{Re
}\lambda_0+\varepsilon}),\label{eq337}
\end{equation}
$\mbox{for all }\varepsilon>0$. Hence (\ref{eq335}) holds as can
be seen from (\ref{eq334}) and (\ref{eq336}) holds by
(\ref{eq337}).

\section{Estimates for Green's function}

We shall consider estimates for the Green's function for the
Lam\'e and Stokes systems on a cone ${\cal K}$ in a spirit of the
three dimensional result of Maz'ya and Plamenevskii \cite{MP} and
Theorem 11.4.7 of \cite{MR}.

Let ${D}\subset {\R}^n$ be a bounded domain that is smooth
everywhere except at a single point (without loss of generality we
can assume this point is $0$). In a small neighborhood of this
point we will assume the domain looks like a cone $\cal K$ defined
above. That is for some $\delta>0$
\begin{equation}
D\cap B(0,\delta)= {\cal K} \cap B(0,\delta),\label{eq13}
\end{equation}
where ${\cal K}=\{x\in {\R}^n; 0<|x|<\infty, x/|x|\in\Omega\}$.
Recall that $\Omega\subset {\mathbb S}^{n-1}$. We will assume that
$\Omega$ is open and nonempty.

Let us denote the vector $\delta_j$ with components
$(\delta_{1j},\delta_{2j},\dots,\delta_{nj})$ where $\delta_{ij}$
denotes the Kronecker symbol. Consider the fundamental solution
$(g^j,p^j)$ of our system in $D$, that is the solution of the
problem
\begin{eqnarray}
-\Delta_x g^j(x,\xi) +\nabla_x
p^j(x,\xi)&=&\delta(x-\xi)\delta_j,\nonumber\\
\mbox{div}_x g^j(x,\xi)+(1-2\nu)p^j(x,\xi)&=&0,\qquad\qquad
x,\xi\in D.\label{eq14}
\end{eqnarray}

We claim that the following holds:

\begin{eqnarray}
|\nabla_{\xi} g^j(x,\xi) |&\le& c|x-\xi|^{-(n-1)},\qquad\mbox{if
}|x|<|\xi|<2|x|,\nonumber\\
|\nabla_{\xi} g^j(x,\xi) |&\le&
c|x|^{\alpha}|\xi|^{-(\alpha+n-1)},\qquad\mbox{if }
2|x|<|\xi|,\label{eq15}\\
|\nabla_{\xi} g^j(x,\xi)
|&\le&c|\xi|^{\alpha-1}|x|^{-(\alpha+n-2)},\qquad\mbox{if }
2|\xi|<|x|.\nonumber
\end{eqnarray}
Here $\alpha>0$ is a small constant as in Theorem
\ref{t3}.\vglue2mm

We claim it suffices to establish estimates (\ref{eq15}) for the
fundamental solution in the unbounded cone ${\cal K}$ with zero
Dirichlet boundary conditions at $\partial{\cal K}$ and infinity.
Here the fundamental solution in the unbounded cone ${\cal K}$ is
a solution of (\ref{eq15}) in ${\cal K}$ that decays at infinity,
i.e., $g^j(x,\xi)\to 0$ as $|x|\to \infty$ sufficiently fast so
that $g^j(.,\xi)$ is an $L^2$ function outside the pole at
$x=\xi$.

We will explain the step of going from ${\cal K}$ to $D$ in detail
below. Let us now work on ${\cal K}$. The existence and uniqueness
of Green's function in an infinite cone for general elliptic
boundary value problems was established in Theorem 7.2 of
\cite{KMR2}, in particular the first estimate (\ref{eq15}) follows
directly from this Theorem.

We now look at the last estimate in (\ref{eq15}). Given the
homogeneity of $g^j$ on ${\cal K}$ we have
$$g^j(x,\xi)=\lambda^{n-2}g^j(\lambda x,\lambda\xi), \mbox{ for all } \lambda>0.$$
It follows that
$$g^j(x,\xi)=|x|^{2-n}g^j(x/|x|,\xi/|x|).$$

Fix now $x\in{\cal K}$. On domain $|\xi|\le 0.99999|x|$,
$\xi\mapsto g^j(x,\xi)$ is just the solution of the adjoint
problem to (\ref{eq01}) - so all we proved for (\ref{eq01}) also
holds for the adjoint equation. In particular (\ref{eq335})
applies and
\begin{eqnarray}
|\nabla_\xi g^j(x,\xi)|&=&|x|^{2-n}|\nabla_\xi
g^j(x/|x|,\xi/|x|)|\le
C|x|^{2-n}\left|\frac{\xi}{x}\right|^{\alpha-1}\frac1{|x|}\nonumber\\
&=&C|\xi|^{\alpha-1}|x|^{-(\alpha+n-2)},\qquad \mbox{for
}|\xi|<3/4|x|.
\end{eqnarray}

The second estimate is similar, but we use (\ref{eq336}). Again we
have
$$g^j(x,\xi)=|x|^{2-n}g^j(x/|x|,\xi/|x|),$$
for on domain domain $|\xi|>1.00001|x|$. As before $\xi\mapsto
g^j(x,\xi)$ solves an adjoint problem, so by (\ref{eq336}) we get
for $|\xi|>4/3|x|$:
\begin{eqnarray}
|\nabla_\xi g^j(x,\xi)|&=&|x|^{2-n}|\nabla_\xi
g^j(x/|x|,\xi/|x|)|\le C
|x|^{2-n}\left|\frac{\xi}{|x|}\right|^{1-n-\alpha}\frac1{|x|}\nonumber\\
&=&C|x|^\alpha|\xi|^{-(\alpha+n-1)}.
\end{eqnarray}

Having required estimates on $\partial {\cal K}$ we show now that
the same will be true on the domain $D$ defined at the beginning
of this section. The point is that by (\ref{eq13}) the domains $D$
and ${\cal K}$ coincide. Hence, if we denote the Green's function
on ${\cal K}$ by $\widetilde{g^j}$, then the Greens's function
$g^j$ for $D$ can be sought in the form
$$g^j(x,\xi)=\widetilde{g^j}(x,\xi)\varphi(x)+f_j(x,\xi),$$
where $\varphi(x)$ is a smooth cut-off function equal to one on
$B(0,\delta/2)$ and vanishing outside $B(0,\delta)$.

Since
\begin{eqnarray}
-\Delta_x [g^j(x,\xi)\varphi(x)] +\nabla_x
p^j(x,\xi)&=&\delta(x-\xi)\delta_j+\rho(x,\xi),\nonumber\\
\mbox{div}_x
[g^j(x,\xi)\varphi(x)]+(1-2\nu)p^j(x,\xi)&=&\tau(x,\xi),\qquad\qquad
x,\xi\in D,\label{eq1444}
\end{eqnarray}
where $\rho$ and $\tau$ are are only supported in
$D\cap\{\delta/2<|x|<\delta\}$, we see that $f_j$ must solve
\begin{eqnarray}
-\Delta_x f_j(x,\xi) +\nabla_x
p^j(x,\xi)&=&-\rho(x,\xi),\nonumber\\
\mbox{div}_x
f_j(x,\xi)+(1-2\nu)p^j(x,\xi)&=&-\tau(x,\xi),\qquad\qquad x,\xi\in
D,\label{eq1445}
\end{eqnarray}
and $f_j\big|_{\partial D}=0$. The main point is that this reduces
the problem to dealing with the remainder term $f_j$. Now however,
near the singular vertex Theorem \ref{t3} applies and away from
the vertex the domain $\partial D$ is smooth, hence $f_j$ is
smooth there as well, so $f_j$ has the required regularity. See
also section 4 of \cite{MP3} for estimates of this type.

\section{The $L^p$ Dirichlet problem}

Let $D$ be the domain defined in the previous section. We would
like to study the solvability of the classical $L^p$ Dirichlet
problem for the Lam\'e and Stokes systems in the domain $D$. Let
us recall the definition.

\begin{definition} Let $1<p\le \infty$. We say that the Dirichlet
problem for the Lam\'e system ($\nu<1/2$) or Stokes system
($\nu=1/2$) is $L^p$ solvable on the domain $D$ if for all vector
fields $f\in L^p(\partial D)$ there is  pair $(u,p)$ such that

\begin{eqnarray}\label{eq16}
-\Delta u+\nabla p&=&0,\qquad \mbox{div
}u+(1-2\nu)p=0\quad\mbox{in
}D\\
u\big|_{\partial D}&=&f\quad\mbox{almost everywhere},\nonumber\\
u^*&\in& L^p(\partial D),\nonumber
\end{eqnarray}

and for some $C>0$ independent of $f$

$$\|u^*\|_{L^p(\partial D)}\le C\|f\|_{L^p(\partial D)}.$$

Here, the boundary values of $u$ are understood in the
nontangential sense, that is we take the limit
$$u\big|_{\partial D}(x)=\lim_{y\to x,\,y\in \Gamma(x)} u(y),$$
over a collection of nontangential cones $\Gamma(x)$ of same
aperture and vertex at $x\in \partial D$ and $u^*$ is the
classical nontangential maximal function defined as
$$u^*(x)=\sup_{y\in \Gamma(x)}|u(y)|,\qquad\mbox{for all }x\in\partial D.$$
\end{definition}

Our main result is

\begin{theorem}\label{t2} Let $D$ be the domain defined above and $\nu\le 1/2$. Then for any
$(n-1)/(\alpha+n-2)<p\le \infty$ ($\alpha$ is same as in
(\ref{eq15})) the $L^p$ Dirichlet problem for the system
(\ref{eq16}) is uniquely solvable. Moreover, for any such $p$
there exists a constant $C(p)>0$ such that the solution $(u,p)$ of
the problem with boundary data $f\in L^p$ satisfies the estimate
$$\|u^*\|_{L^p(\partial D)}\le C(p) \|f\|_{L^p(\partial D)}.$$
Moreover, if $f\in C(\partial D)$, then $u\in
C(\overline{\Omega})$ and an estimate
$$\|u\|_{C(\overline{D})}\le C\|f\|_{C(\partial D)}$$
holds.
\end{theorem}

\noindent{\it Proof:} Consider the representation of the solution
$u$ by the Green's formula. That is, for $j=1,2,\dots,n$
\begin{equation}
u_j(x)=-\int_{\partial D}\frac{\partial
g^j(x,\xi)}{\partial\nu_\xi}f(\xi)d\sigma_\xi,\qquad x\in D,
\end{equation}
where $\nu_\xi$ is the outer normal at the boundary point $\xi$
and $d\sigma$ is the $(n-1)$-dimensional area element at $\partial
D$. We will use (\ref{eq15}) to establish the claim. From this
estimate we obtain in the zone $|x|/2<|\xi|<2|x|$
$$\left|\frac{\partial
g^j(x,\xi)}{\partial\nu_\xi}\right|\le C\frac{R(x)}{|x-\xi|^n},$$
where $R(x)=\text{dist}(x,\partial D)$. From this and other two
estimates of (\ref{eq15}) we obtain:
\begin{eqnarray}\nonumber
|u_j(x)|&\le&
C\left(|x|^\alpha\int_{E_1}\frac{|f(\xi)|}{|\xi|^{n-1+\alpha}}d\sigma+R(x)\int_{E_2}\frac{|f(\xi)|}{|x-\xi|^n}d\sigma\right.
\\&+&\left.\frac1{|x|^{n-2+\alpha}}\int_{E_3}\frac{|f(\xi)|}{|\xi|^{1-\alpha}}d\sigma\right).\label{eq17}
\end{eqnarray}
Let us denote these three integrals by $v^1(x)$, $v^2(x)$,
$v^3(x)$, respectively. Here
\begin{eqnarray}\nonumber
E_1&=&\{\xi\in\partial D:2|x|<|\xi|\},\\
E_2&=&\{\xi\in\partial D:|x|/2\le|\xi|\le 2|x|\}\label{eq18}\\
E_3&=&\{\xi\in\partial D:2|\xi|<|x|\}.\nonumber
\end{eqnarray}
We deal with these three terms separately. We introduce
$$v^{i,*}(x)=\sup_{y\in \Gamma(x)}|v^i(y)|,\qquad\mbox{for all }x\in\partial D\mbox{ and }i=1,2,3.$$
It follows that
$$u^*(x)\le C(v^{1,*}(x)+v^{2,*}(x)+v^{3,*}(x))\qquad\mbox{for all }x\in\partial D.$$

\begin{lemma}\label{l1} There exists $C>0$ such that
\begin{eqnarray}\label{eq19}
\|v^{1,*}\|_{L^\infty(\partial D)}&=&\|v^{1}\|_{L^\infty(D)}\le C
\|f\|_{L^\infty(\partial D)}\\
\|v^{1,*}\|_{L^{1,w}(\partial D)}&\le& C\|f\|_{L^1(\partial
D)}.\nonumber
\end{eqnarray}
Here $L^{p,w}$ is the weak-$L^p$ space equipped with the norm
$$\|f\|_{L^{p,w}}=\left(\sup_{\lambda>0}\lambda^p\sigma(\{\xi:|f(\xi)|>\lambda\})\right)^{1/p}.$$
\end{lemma}

\noindent{\it Proof of the lemma:} The definition of $v^1$ implies
that
\begin{eqnarray}
v^1(x)&\le&
|x|^\alpha\int_{E_1}\frac{|f(\xi)|}{|\xi|^{n-1+\alpha}}d\sigma\nonumber\\&\le&
|x|^\alpha\|f\|_\infty\int_{|\xi|>2|x|}|\xi|^{-(n-1+\alpha)}d\sigma\le
C\|f\|_\infty,\label{eq20}
\end{eqnarray}
since the integral is bounded by $C|x|^{-\alpha}$. From this the
first claim follows. On the other hand, when $f\in L^1$ we use a
trivial estimate $|\xi|^{-(n-1+\alpha)}\le C|x|^{-(n-1+\alpha)}$
which gives us
\begin{equation}
v^1(x)=
|x|^\alpha\int_{E_1}\frac{|f(\xi)|}{|\xi|^{n-1+\alpha}}d\sigma\le
C|x|^{-(n-1)}\int_{E_1}|f(\xi)|\sigma\le
C\|f\|_1|x|^{-(n-1)}.\label{eq21}
\end{equation}
From this we can estimate $v^{1,*}(\xi)$ for $\xi\in \partial D$.
We realize that for any $x\in \Gamma(\xi)$ we have $|x|\ge |\xi|$,
hence
$$v^{1,*}(\xi)\le C\|f\|_1|\xi|^{-(n-1)}.$$
It follows that
\begin{eqnarray}\label{eq22}
\sigma(\{\xi: v^{1,*}(\xi)>\lambda\})&\le&
\sigma(\{\xi:C\|f\|_{L^1}|\xi|^{-(n-1)}>\lambda\})\\&=&\nonumber\sigma\left(\left\{\xi:|\xi|^{n-1}<\frac{C\|f\|_{L^1}}{\lambda}\right\}\right)\le
C\frac{\|f\|_{L^1}}{\lambda},
\end{eqnarray}
since the surface measure of a ball $\{\xi:|\xi|\le R\}$ is
proportional to $R^{n-1}$. From this the fact that $v^{1,*}$
belongs to the weak-$L^1$ follows.

\begin{lemma}\label{l2} There exist $C(p)>0$ such that for all $p>\frac{n-1}{n-2+\alpha}$
\begin{eqnarray}\label{eq23}
\|v^{3,*}\|_{L^\infty(\partial D)}&=&\|v^{3}\|_{L^\infty(D)}\le C
\|f\|_{L^\infty(\partial D)}\\
\|v^{3,*}\|_{L^{p,w}(\partial D)}&\le& C\|f\|_{L^p(\partial
D)}.\nonumber
\end{eqnarray}
\end{lemma}

\noindent{\it Proof of the lemma:} The definition of $v^3$ implies
that for $q=p/(p-1)$
\begin{eqnarray}
v^3(x)&=&
|x|^{-(n-2+\alpha)}\int_{E_3}\frac{|f(\xi)|}{|\xi|^{1-\alpha}}d\sigma\nonumber\\&\le&
C|x|^{-(n-2+\alpha)}\|f\|_{L^p(E_3)}\left(\int_{E_3}\frac1{|\xi|^{q(1-\alpha)}}d\sigma\right)^{1/q}.
\label{eq24}
\end{eqnarray}
In polar coordinates
\begin{equation}
\int_{E_3}\frac1{|\xi|^{q(1-\alpha)}}d\sigma\approx\int_0^{|x|/2}r^{n-2-q(1-\alpha)}dr<\infty
\label{eq25}
\end{equation}
if and only if $n-2-q(1-\alpha)>-1$ or $p>\frac{n-1}{n-2+\alpha}$.
Assuming that (\ref{eq24}) gives us that
\begin{equation}
v^3(x)\le
C|x|^{-(n-2+\alpha)}\|f\|_{L^p(E_3)}\left(|x|^{n-1-q(1-\alpha)}\right)^{1/q}\le
C\|f\|_{L^p(\partial D)}|x|^{-\frac{n-1}p}. \label{eq26}
\end{equation}
When $p=\infty$ the first part of (\ref{eq23}) follows. For
$\frac{n-1}{n-2+\alpha}<p<\infty$ we observe that as before
$$v^{3,*}(\xi)=\sup_{x\in \Gamma(\xi)}v^3(x)\le C\|f\|_{L^p(\partial D)}|\xi|^{-\frac{n-1}p}.$$
Hence
\begin{eqnarray}\label{eq27}
\sigma(\{\xi: v^{3,*}(\xi)>\lambda\})&\le&
\sigma(\{\xi:C\|f\|_{L^p}|\xi|^{-(n-1)/p}>\lambda\})\\&=&\nonumber\sigma\left(\left\{\xi:|\xi|^{n-1}<\frac{C\|f\|^p_{L^p}}{\lambda^p}\right\}\right)\le
C\frac{\|f\|^p_{L^p}}{\lambda^p}.
\end{eqnarray}
This gives the second estimate of (\ref{eq23}).

\begin{lemma}\label{l3} There exist $C>0$ such that
\begin{eqnarray}\label{eq28}
\|v^{2,*}\|_{L^\infty(\partial D)}&=&\|v^{2}\|_{L^\infty(D)}\le C
\|f\|_{L^\infty(\partial D)}.
\end{eqnarray}
\end{lemma}
\noindent{\it Proof of the lemma:} For any $x\in D$:
\begin{equation}\label{eq29}
v^2(x)=R(x)\int_{E_2}\frac{|f(\xi)|}{|x-\xi|^n}d\sigma\le
\|f\|_{L^\infty(D)} R(x)\int_{E_2}|x-\xi|^{-n}d\sigma.
\end{equation}
We need to consider how the set $E_2$ looks. Clearly for every
$x\in E_2$, $|x-\xi|\ge R(x)$. Since $E_2=\{\xi\in\partial
D:|\xi|\in [|x|/2,2|x|]\}$ we can parameterize $E_2$ and think
about it as a cylinder $B\times [0,A]$, where $B$ is an
$n-2$-dimensional set of diameter at most $2|x|$. In this
parametrization for $\xi=(b,s)\in B\times [0,A]$ we have $|x-\xi
|\approx R(x)+s$. It follows that
\begin{eqnarray}\label{eq30}
&&R(x)\int_{E_2}|x-\xi|^{-n}d\sigma \approx
R(x)\int_0^A\frac{R(x)^{n-2}}{(R(x)+s)^n}ds\\&\le&
R(x)^{n-1}\int_{R(x)}^\infty s^{-n}ds\le C.\nonumber
\end{eqnarray}
Hence that
$$v^2(x)\le C\|f\|_{L^\infty(D)},$$
with constant $C>0$ independent of the point $x\in D$.\vglue1mm

To handle $p<\infty$ we need to introduce further splitting.
Recall that $v^{2,*}(x)$ for a boundary point $x\in\partial D$ is
defined as a supremum of $v^2$ over a nontangential cone
$\Gamma(x)$ with vertex at $x$. The points $y\in \Gamma(x)$ are of
two kinds. The first kind are points for which $|y-x|\le |x|$
(these are near the vertex $x$). The second kind are points
$|y-x|>|x|$, for these we can make a simple observation that
$R(y)\approx |y-x|\approx |y|$. To distinguish these two we
introduce
\begin{eqnarray}\label{eq31}
w(x)&=&\sup\{v^2(y): y\in\Gamma(x)\mbox{ and }|y-x|\le |x|\},\\
z(x)&=&\sup\{v^2(y): y\in\Gamma(x)\mbox{ and
}|y-x|>|x|\}.\nonumber
\end{eqnarray}
It follows that pointwise $v^{2,*}(x)\le w(x)+z(x)$ for
$x\in\partial D$ and hence
$$\|v^{2,*}\|_{L^p(\partial D)}\le \|w\|_{L^p(\partial D)}+\|z\|_{L^p(\partial D)},\qquad\mbox{for any }1\le p\le\infty.$$

Let us denote by $B_{a,b}$ the part of the boundary of $\partial
D$ such that
\begin{equation}\label{eq32}
B_{a,b}=\{\xi\in\partial D; a\le |\xi|\le b\}\qquad\mbox{for
}0<a<b.
\end{equation}
Due to our assumption that near $0$ our domain looks like a cone,
it follows that for $a,b,\lambda b<\delta$ the sets $B_{a,b}$ and
$B_{\lambda a,\lambda b}$ ($\lambda>0$) are simple rescales of
each other, that is
\begin{equation}\label{eq33}
B_{\lambda a,\lambda b}=\lambda B_{a,b},
\end{equation}
where the multiplication of a set by a scalar is understood in the
usual sense. We claim that for any $\lambda>0$ and $1<p\le\infty$
we have an estimate
\begin{equation}\label{eq34}
\|w\|_{L^p(B_{\lambda,2\lambda})}\le C_p
\|f\|_{L^p(B_{\lambda/2,8\lambda})}.
\end{equation}
It is enough to establish this for a single value of $\lambda>0$,
since then due to the rescaling argument the statement must be
true for all $\lambda>0$ small as follows from (\ref{eq33}). The
(\ref{eq34}) holds, since the Dirichlet problem for the Lam\'e (or
Stokes) system is solvable for all $1<p\le\infty$, provided the
domain is $C^1$. As the sets $B_{\lambda,2\lambda}$ and
$B_{\lambda/2,8\lambda}$ are outside the singularity (vertex at
$0$), our domain can be modified near the vertex outside these
sets, so that the whole domain is $C^1$. Then the solvability for
all $p>1$ is used to get (\ref{eq34}).

Setting $\lambda=2^{-n}\delta$ and summing over $n$ we get:
\begin{eqnarray}\label{eq35}
\|w\|^p_{L^p(\partial D\cap B(0,\delta/8))}&=&\sum_{n=8}^\infty
\|w\|^p_{L^p(B_{2^{-n-1}\delta,2^{-n}\delta})}\\&\le& C_p^p
\sum_{n=8}^\infty
\|f\|^p_{L^p(B_{2^{-n-2}\delta,2^{-n+3}\delta})}\le
4C_p^p\|f\|^p_{L^p(\partial D\cap B(0,\delta))}.\nonumber
\end{eqnarray}
This is the necessary estimate for $w$.

Looking at $z$, let us pick a point $y\in \Gamma(x)$ for
$x\in\partial D$ such that $|y-x|>|x|$. We need to estimate
$v^2(y)$:
$$v^2(y)=R(y)\int_{E_2}\frac{|f(\xi)|}{|y-\xi|^n}d\sigma.$$
Clearly, due to the fact that $R(y)\approx |y|$ and since for
$\xi\in E_2$: $|\xi|\approx |y|$ we see that for any $\xi\in E_2$
we have $|y-\xi|\approx |y|$. Hence
$$v^2(y)\le C|y|^{1-n}\int_{E_2}|f(\xi)|d\sigma\le
C|y|^{1-n}\|f\|_{L^1(\partial D)}.$$ It follows that for
$x\in\partial D$:
\begin{equation}
z(x)=\sup\{v^2(y): y\in\Gamma(x)\mbox{ and }|y-x|>|x|\}\le
C|x|^{1-n}\|f\|_{L^1(\partial D)}.\label{eq333}
\end{equation}
 By the same argument as in
(\ref{eq22}) it follows that $z$ belongs to a weak
$L^{1,w}(\partial D)$. Hence we can claim that

\begin{lemma}\label{l4} For any $1<p\le\infty$ there exists $C_p>0$
such that
\begin{eqnarray}\label{eq36}
\|v^{2,*}\|_{L^p(\partial D)}\le C_p \|f\|_{L^p(\partial D)}.
\end{eqnarray}
\end{lemma}

\noindent{\it Proof. } Consider first a mapping $f\mapsto z_f$,
where for given $f$, we define $z=z_f$ by (\ref{eq333}). This is
not a linear mapping, but it is sublinear, that is
$$z_{f+g}\le z_f+z_g.$$
By Lemma \ref{l3} this mapping is bounded on $L^\infty$ and also
as we have just show maps $L^1$ to weak $L^{1,w}$. By the
Marcinkiewicz interpolation theorem (which only requires
sublinearity, not linearity) this mapping is therefore bounded on
any $L^p$, $p>1$, and we have the estimate:
$$\|z_f\|_{L^p}\le C_p\|f\|_{L^p}.$$
As we already know this for $w$ we see that
$$\|v^{2,*}\|_{L^p}\le \|w\|_{L^p}+\|z\|_{L^p}\le C_p\|f\|_{L^p},$$
for all $p>1$.\vglue1mm

\noindent{\it Proof of Theorem \ref{t2}.} We use the Marcinkiewicz
interpolation theorem in the same spirit as we did above. By
Lemmas \ref{l1} and \ref{l2} it follows that for all $p>1$:
$$\|v^{1,*}\|_{L^p}\le C_p\|f\|_{L^p},$$
and for all $p>\frac{n-1}{n-2+\alpha}$:
$$\|v^{3,*}\|_{L^p}\le C_p\|f\|_{L^p}.$$
Combining the estimates for $v^{1,*}$, $v^{2,*}$ and $v^{3,*}$
yields the desired claim, since
$$\|u^*\|_{L^p}\le C(\|v^{1,*}\|_{L^p}+\|v^{2,*}\|_{L^p}+\|v^{3,*}\|_{L^p}).$$

\end{document}